\documentclass[12pt]{article} 

\usepackage{hyperref}
\usepackage{amsmath, amsthm, amssymb, algorithm, algorithmic,
  enumerate, verbatim}

\theoremstyle{plain}
\newtheorem{theorem}{Theorem}[section]
\newtheorem{proposition}[theorem]{Proposition}
\newtheorem{lemma}[theorem]{Lemma}

\theoremstyle{remark} 
\newtheorem{example}[theorem]{Example}

\newcommand{\DMO}{\DeclareMathOperator}
\newcommand{\fl}[2]{F_{< #2}(#1)}
\DMO{\lspec}{lspec}
\DMO{\lcm}{lcm}

\title{Length Spectra of Natural Numbers} 
\author{Wai Yan Pong\\ California State University Dominguez Hills}

\begin{document}
\maketitle
\begin{abstract} 
  Two numbers are spectral equivalent if they have the same length
  spectrum. We show how to compute the equivalence classes of this
  relation. Moreover, we show that these classes can only have either
  $1,2$ or infinitely many elements.
\end{abstract}

\section{Introduction} 
\label{s:intro} Some numbers can be written as a sum of consecutive
integers, for example, $9=2+3+4$ but some cannot, for example, $8$.
The following is a beautiful characterization of this phenomenon:
\begin{theorem} 
  \label{th:po2} A number is a sum of consecutive integers if and only
  if it is not a power of 2.
\end{theorem} 
Proofs of Theorem~\ref{th:po2} can be found in~\cite{guy, lev,
  sumcon}.

To simplify the subsequent discussion, let us call a sequence of
consecutive natural numbers a {\bf decomposition of $\pmb n$} if its
terms sum to $n$. The {\bf length} of a decomposition is the number of
terms in the decomposition and the {\bf parity} of a decomposition is
the parity of its length. A {\bf trivial} decomposition is a
decomposition of length $1$. Clearly every number $n$ has a trivial
decomposition, namely $(n)$. The following result in~\cite{sumcon}
(also in~\cite{guy, lev}) is fundamental and will be used frequently
throughout this article.
\begin{theorem} 
  \label{th:corr} Let $n$ be a natural number and $k$ be an odd factor
  of $n$. If $k^2 < 2n$, then the sequence
  \begin{equation} 
    \label{eq:odd} \frac{n}{k}-\frac{k-1}{2},\quad
    \frac{n}{k}-\frac{k-1}{2}+1, \quad \cdots, \quad
    \frac{n}{k}+\frac{k-1}{2}
  \end{equation} 
  is an odd decomposition of $n$ of length $k$. On the other hand, if
  $k^2 > 2n$, then the sequence
  \begin{equation} 
    \label{eq:even} \frac{k-1}{2}- \frac{n}{k} +1,
    \quad \frac{k-1}{2}-\frac{n}{k} + 2, \quad \cdots, \quad \frac{n}{k} +
    \frac{k-1}{2}
  \end{equation} 
  is an even decomposition of $n$ of length $2n/k$. Moreover, every
  decomposition of $n$ has one of these forms.
\end{theorem} 
Theorem~\ref{th:corr} explicitly demonstrates a 1-to-1 correspondence
between the odd factors of $n$ and its decompositions.  Since the
powers of $2$ are the only numbers having no odd factors other than
$1$, their decompositions can only be trivial. This establishes
Theorem~\ref{th:po2} as a consequence of Theorem~\ref{th:corr}.

The {\bf length spectrum} (or simply the {\bf spectrum}) of a number
$n$, denoted by $\lspec(n)$, is the set of lengths of the
decompositions of $n$. According to Theorem~\ref{th:corr}, the
spectrum of $n$ is the set
\begin{equation} 
  \label{eq:lspec} \{k \colon k\ \text{odd},\ k \mid n,\ k^2 < 2n\} 
  \cup \{2n/k \colon k\ \text{odd},\ k \mid n,\ k^2 > 2n\}.
\end{equation} 
As an example, we list in Table~\ref{tb:lspec(45)} the decompositions
of the number $45$ along with their lengths, parities and associated
odd factors.
\begin{table}[htp]
\begin{center}
\begin{tabular}{|c|c|c|c|} 
  \hline 
  factor & decomposition & length & parity\\ 
  \hline 
  1 &(45) &1 &odd\\ 45 &(22, 23) &2 &even\\ 
  3 &(14, 15, 16) &3 &odd\\ 
  5 &(7, 8, 9, 10, 11) &5 &odd\\ 
  15 &(5, 6, 7, 8, 9, 10) &6 &even\\ 
  9 &(1, 2, 3, 4, 5, 6, 7, 8, 9) &9 &odd\\ 
  \hline
\end{tabular}
\caption{The spectrum of $45$ is $\{1,2,3,5,6,9\}$}
\label{tb:lspec(45)}
\end{center}
\end{table}

We say that two numbers are {\bf spectral equivalent} if they have the
same length spectrum. The {\bf spectral class} of $n$, denoted by
$L(n)$, is the equivalence class of $n$ under spectral equivalence.
With these notions, Theorem~\ref{th:po2} can be restated in the
following way: The powers of 2 form a spectral class with $\{1\}$ as
their common spectrum. The odd primes form another interesting
spectral class. According to Theorem~\ref{th:corr}, they are the
numbers having $\{1,2\}$ as spectrum. These two examples motivate the
following question:
\begin{quote} 
  \em Given $n$, can we compute its spectral class?
\end{quote} 
We will give an algorithm in Section~\ref{s:e-a} which answers this
question in an affirmative way. We also derive from it another
algorithm which solve the following problem.
\begin{quote} 
  \em Given a finite set of numbers $S$, compute the set of numbers
  with $S$ as their common spectrum.
\end{quote} 

\section{Properties of Length Spectra}
\label{s:lspec} The following notations and conventions will be
adopted throughout the rest of this article.
\begin{itemize}
\item For a set $A$, we write $|A|$ for its cardinality.
  
\item For a set of natural numbers $A$, we write $A_0$ and $A_1$ for
  the set of even and odd elements of $A$, respectively.

\item For a rational number $c$ and a set of natural numbers $A$, we
write $cA$ for the set $\{ ca \colon a \in A\}$.
  
\item For a prime $p$ and a natural number $m$, we write $v_p(m)$ for
  the exponent of $p$ in the prime factorization of $m$.
\end{itemize}
  
Let $n$ be a fixed but arbitrary natural number, we use
$(k_i)_{i=1}^s$ to denote the list of odd factors of $n$ in ascending
order.  So $k_1=1$ and $n=2^{\nu}k_s$ where $\nu=v_2(n)$. We use $r$
to denote the largest index such that $k_r^2 < 2n$. Let us note that
for $1\le i \le s$, $k_j$ and $k_{s-j+1}$ are complementary factors of
$k_s$, i.e.  $k_jk_{s-j+1} = k_s$ hence the length spectrum of $n$ can
be re-written as
\begin{equation}
  \label{eq:newform} \{k_i \colon 1 \le i \le r\}
  \cup 2^{\nu+1}\{k_i \colon 1 \le i \le s-r\}.
\end{equation}

The first important observation about length spectra is:
\begin{theorem}
  \label{th:even<=odd} The number of even decompositions of a number
  is at most the number of its odd decompositions.
\end{theorem}
\begin{proof} 
  Suppose $k_j$ corresponds to an even decomposition of $n$. By
  Theorem~\ref{th:corr}, $k_j^2 > 2n$ and so $k_{s-j+1}^2 =
  k_s^2/k_j^2 < (2n)^2/2n =2n$. Therefore, again by
  Theorem~\ref{th:corr}, $k_{s-j+1}$ corresponds to an odd
  decomposition of $n$. The theorem now follows since the map $k
  \mapsto k_s/k$ is 1-to-1.
\end{proof} 
In the light of Theorem~\ref{th:even<=odd}, we call a spectrum
\begin{itemize}
\item {\bf balanced} if it has an equal number of even and odd
  elements;

\item {\bf unmixed} if it has no even elements;

\item {\bf lopsided} if it is neither balanced nor unmixed.
\end{itemize} 
Our next result characterizes the numbers with an unmixed spectrum.
\begin{theorem}
  \label{th:class-unmixed} The set of numbers with an unmixed spectrum
  is
\[ 
\{ 2^{\alpha}k \colon \alpha \ge 0,\ k \text{ odd, } 2^{\alpha+1} > k
\}.
\]
\end{theorem}
\begin{proof} 
  It follows from~\eqref{eq:lspec} that a number with no even
  decompositions if and only if it is of the form $2^{\alpha}k$ with
  $\alpha \ge 0$, $k$ odd and $2^{\alpha+1}k > k^2$, i.e.
  $2^{\alpha+1} > k$.
\end{proof}

Numbers with a balanced spectrum are trickier to capture. For a
natural number $k$, let $q(k)$ be the minimum of the ratios $m'/m$
where $m \le m'$ are complementary factors of $k$. Note that by
definition $q(k) \ge 1$. Moreover, for a pair $m,m'$ of complementary
factors of $k$, $m'/m=q(k)$ if and only if no factor of $k$ is
strictly in between $m$ and $m'$. Also, it may be worth pointing out
that $q(k)$ measures how far is $k$ from being a prefect square: $k$
is a prefect square if and only if $q(k)=1$.
\begin{theorem} 
  \label{th:class-balanced} The set of numbers with a balanced
  spectrum is
  \[ 
  \{2^{\alpha}k \colon \alpha \ge 0,\ k \text{ odd, } q(k) >
  2^{\alpha+1} \}.
  \]
\end{theorem}
\begin{proof} 
  Suppose $2^{\alpha}k$ is a member of the set in display.  Let $m,
  m'$ be a pair of complementary factors of $k$ with $m'/m=q(k)$. Then
  $m'/m > 2^{\alpha+1}$, hence $m' > 2^{\alpha+1}m$ and so
  \begin{equation}
    \label{eq:ineq} 
    m'^2 > 2^{\alpha+1}mm' = 2^{\alpha+1}k > m^2.
  \end{equation} 
  Since no factor of $k$ is strictly in between $m$ and $m'$, the
  inequalities in~\eqref{eq:ineq} imply the odd (resp. even)
  decompositions of $2^{\alpha}k$ correspond precisely to those
  factors of $k$ that are $\le m$ (resp. $\ge m'$).  Thus the
  assignment $l \mapsto mm'/l$ induces an injective map from the odd
  decompositions of $2^{\alpha}k=2^{\alpha}mm'$ to its even
  decompositions. Therefore, by Theorem~\ref{th:even<=odd},
  $2^{\alpha}k$ has a balanced spectrum.

  Conversely, suppose $n$ has a balanced spectrum. Then we have (in
  the notations introduced earlier) $s=2r$. Hence
  $n=2^{\nu}k_rk_{r+1}$ and $k_{r+1}/k_r = q(k_rk_{r+1})$.  Moreover,
  $k_{r+1}^2 > 2n= 2^{\nu+1}k_rk_{r+1}$, therefore $k_{r+1}/k_r >
  2^{\nu+1}$ and so $n$ belongs to the set displayed in the statement
  of the theorem.
\end{proof}

\section{Spectral Classes} 
\label{s:lsc} In this section, we determine the spectral class of a
number according to the type of its spectrum.  We begin with an
observation which is clear from the form of the spectrum given
in~\eqref{eq:newform}:
\begin{proposition} 
  \label{p:even} If $n$ has an even decomposition, then the highest
  power of $2$ dividing $n$ is half the least even element of
  $\lspec(n)$. In fact, $v_2(n)=v_2(e)-1$ for any even element $e$ of
  $\lspec(n)$.
\end{proposition}

Let $S$ be a finite set of natural numbers. Recall that $S_i$ is the
set of elements of $S$ congruent to $i \pmod{2}$. We define $D(S)$,
the {\bf difference set} of $S$, to be $S_1$ with its least $|S_0|$
elements removed if $|S_0| \le |S_1|$, or the empty set otherwise.
\begin{proposition}
  \label{p:largest} If a number has more odd than even decompositions,
  then its greatest odd factor is the product of the maximum and
  minimum of the difference set of its spectrum.
\end{proposition}
\begin{proof} 
  Suppose $S=\lspec(n)$ and $|S_1| > |S_0|$. In this case, $D(S)$ is
  the non-empty set $S_1 \setminus 2^{-(\nu+1)}S_0 = \{k_{s-r+1},
  \ldots, k_r\}$. Since $k_{s-r+1}k_r =k_s$, the largest odd factor of
  $n$, the proposition follows.
\end{proof}

\begin{theorem}
  \label{th:lsc-lopsided} If $lspec(n)$ is lopsided, then $L(n) =
  \{n\}$.
\end{theorem}
\begin{proof} 
  If $\lspec(n)$ is lopsided, then by Proposition~\ref{p:even}
  and~\ref{p:largest} both $v_2(n)$ and $k_s$, and therefore $n$, can
  be recovered from $\lspec(n)$.
\end{proof} 
To illustrate the results that we have just proved, let us decide
whether the set $S:=\{1,3,4,5,9,12\} = \{1,3,5,9\} \cup \{4,12\}$ is a
spectrum. First, let us note that $D(S) = \{5,9\}$. So if $S$ is a
spectrum, say $S=\lspec(n)$, then by Proposition~\ref{p:largest} the
largest odd factor of $n$ is $5 \cdot 9 = 45$. By
Proposition~\ref{p:even}, $v_2(n)=v_2(4)-1 =1$.  Therefore, $n$ can
only be $90$. One then finds that $S$ is indeed a spectrum by
verifying $\lspec(90)=S$.

Next we determine the spectral class of $n$ when its spectrum is
unmixed.
\begin{theorem} 
  \label{th:lsc-unmixed} If $\lspec(n)$ is unmixed, then
  \[ 
  L(n)=\{2^{\alpha}m \colon \alpha \ge 0,\ 2^{\alpha+1} > m\} 
  \]
  where $m$ is the largest odd element of $\lspec(n)$.
\end{theorem}
\begin{proof} 
  Suppose $\lspec(n)$ is unmixed. Let $m$ be the largest odd element
  of $\lspec(n)$, $S$ be the set of factors of $m$ and $R$ be the set
  on the right-hand-side of the equation. First, since $\lspec(n)$ is
  unmixed, it follows easily from~\eqref{eq:newform} that
  $\lspec(n)=S$. By Theorem~\ref{th:corr}, every element of $R$ has
  spectrum $S$. By Theorem~\ref{th:class-unmixed}, the converse is
  true hence $R$ is precisely the set of number with spectrum $S$ and
  the theorem follows.
\end{proof}

For a finite set of natural numbers $S$, we define the {\bf
  exceptional set} of $S$ to be the set
\[ 
  E(S)=\{a \in S_1^2 \colon a > m_0,\ \fl{am_1}{a} = S_1\} 
\] 
where $S_1^2$ is the set $\{bc \colon b, c \in S_1\}$, $m_i = \max
S_i$ ($i=0,1$) and $\fl{l}{k}$ denotes the set of factors of $l$ which
are strictly less than $k$. The following simple facts about elements
of exceptional sets will come in handy for our subsequent arguments:
\begin{lemma}
  \label{l:facts} Suppose $S$ is a balanced spectrum and $a \in E(S)$,
  then
  \begin{enumerate}[(i)]
  \item \label{i:notprime} $a$ has a proper prime factor; in
    particular, $a$ is not a prime.
  \item \label{i:proper} every proper factor of $a$ is in $S_1$.
  \item \label{i:Q} $a/m_1 = q(m_1a)$.
  \end{enumerate}
\end{lemma}
\begin{proof} 
  By definition, $a > m_0 \ge 2m_1 \ge 2$. Moreover, $a \in S_1^2
  \setminus S_1$ therefore $a$ cannot be a prime and
  so~\eqref{i:notprime} follows.  Every proper factor of $a$ is
  clearly a member of $\fl{m_1a}{a}$ which is $S_1$ and
  so~\eqref{i:proper} follows.  To show~\eqref{i:Q}, it suffices to
  show that no factor of $m_1a$ is strictly in between $m_1$ and $a$.
  If not, then $\fl{m_1a}{a}$, i.e. $S_1$ will contain a number larger
  than $m_1$, a contradiction.
\end{proof}

The next lemma is crucial to our analysis of balanced spectra.
\begin{lemma}
  \label{l:incl} Suppose $S=\lspec(n)$ is balanced, then
  \[ 
  \frac{m_0}{2}E(S) \subseteq L(n) \subseteq \frac{m_0}{2}(P(S) \cup
  E(S)).
  \] 
  where $P(S)$ is the set of primes that are larger than $m_0$.
\end{lemma}
\begin{proof} 
  Since $S=\lspec(n)$ is balanced, it is of the form $S_1 \cup
  2^{\nu+1}S_1$ where $\nu=v_2(n)$. In particular, $m_0 =
  2^{\nu+1}m_1$.

  Pick $a \in E(S)$, and let $n'$ be $m_0a/2 = 2^{\nu}m_1a$. By
  Lemma~\ref{l:facts}~\eqref{i:Q}, $a/m_1 = q(m_1a)$ and so the
  inequality $a > m_0 = 2^{\nu+1}m_1$ implies $q(m_1a) > 2^{\nu+1}$.
  Therefore, $\lspec(n')$ is balanced according to
  Theorem~\ref{th:class-balanced}. Note that the inequality $a > m_0$
  also implies
  \begin{equation*} 
    a^2 > m_0a = 2n' = 2^{\nu+1}m_1a > m_1^2.
  \end{equation*} 
  These inequalities together with Theorem~\ref{th:corr} and the fact
  that $a/m_1= q(m_1a)$ imply the set of odd elements of $\lspec(n')$
  is $\fl{m_1a}{a}$. But $\fl{m_1a}{a}$ is simply $S_1$, as $a$ is a
  member of $E(S)$. Therefore,
  \[ 
  \lspec(n') = \lspec(n')_1 \cup 2^{v_2(n')+1}\lspec(n')_1 = S_1 \cup
  2^{\nu+1}S_1 = \lspec(n).
  \] 
  So we conclude that $n' \in L(n)$.

  To show the second inclusion, suppose $n' \in L(n)$, i.e.
  $\lspec(n')=S_1 \cup 2^{\nu+1}S_1$. Thus $m_1$ is the largest odd
  factor of $n'$ corresponding to an odd decomposition. Since
  $\lspec(n')$ is balanced, the largest odd factor of $n'$ is of the
  form $m_1a$ for some odd number $a > m_1$ such that $a/m_1 =
  q(m_1a)$. By Proposition~\ref{p:even}, $v_2(n')=\nu$ and so $n' =
  2^{\nu}m_1a=m_0a/2$. By Theorem~\ref{th:class-balanced}, $q(m_1a) >
  2^{\nu+1}$. Therefore, $a > 2^{\nu+1}m_1 = m_0$. Since no factor of
  $m_1a$ is strictly in between $a$ and $m_1$, $\fl{m_1a}{a}$ is the
  set of odd factors of $n'$ not exceeding $m_1$ and that is $S_1$.
  Moreover, since every proper factor of $a$ belongs to $\fl{m_1a}{a}
  =S_1$, that means either $a \in S_1^2$ or $a$ is a prime.  In the
  former case, $a \in E(S)$ and in the latter case, $a \in P(S)$.
\end{proof}

Before stating our next result, which gives the spectral class of $n$
when its spectrum is balanced, we need to introduce one more concept.
Suppose $S$ is a spectrum. By Theorem~\ref{th:corr}, $S_1$ contains
the factors of its elements. In particular, the sets of factors of
$m_1$ is always a subset of $S_1$. In the light of this observation,
we call a balanced spectrum $S$ {\bf non-excessive} if $S_1$ is
precisely the set of factors of $m_1$; otherwise we call $S$ {\bf
  excessive}. Note that we use the word (non-)excessive to describe
balanced spectra only.
\begin{theorem} 
  \label{th:lsc-balanced} Suppose $S=\lspec(n)$ is balanced. Then
  \begin{enumerate}[(i)]
  \item \label{i:non-ex} $L(n)=\dfrac{m_0}{2}(P(S) \cup E(S))$, if $S$
    is non-excessive; or
  \item \label{i:ex} $L(n) = \dfrac{m_0}{2}E(S)$, if $S$ is excessive.
  \end{enumerate}
\end{theorem}
\begin{proof} 
  To proof~\eqref{i:non-ex}, thanks to Lemma~\ref{l:incl}, we only
  need to show that the set $m_0P(S)/2$ is a subset of $L(n)$. Take $p
  \in P(S)$ and let $n'=m_0p/2=2^{\nu}m_1p$ where $\nu=v_2(n)$.  Since
  $p$ is prime, the odd factors of $n'$ are the factors of $m_1$ and
  their $p$ multiples. The inequalities $p > m_0=2^{\nu+1}m_1 > m_1$
  ensure 
  \[
    p^2 > m_0p = 2n' > m_1^2.
  \]
  Therefore, the factors of $m_1$ correspond to the odd decompositions
  of $n'$ while their $p$ multiplies correspond to the even
  decompositions of $n'$ (Theorem~\ref{th:corr}). Thus, $\lspec(n')$
  is balanced and the set of odd elements of $\lspec(n')$ coincides
  with the set of factors of $m_1$.  Since $S$ is non-excessive, the
  set of factors of $m_1$ equals $S_1$.  So we actually have
  \[ \lspec(n') = \lspec(n')_1 \cup 2^{v_2(n')+1}\lspec(n')_1 = S_1
  \cup 2^{\nu+1}S_1 =\lspec(n).
  \] 
  This finishes the proof of Part~\eqref{i:non-ex}.  Incidentally, the
  argument above also shows that if $S=\lspec(n)$ is balanced then
  every element of $m_0P(S)/2$ has non-excessive spectrum.
  Consequently, $L(n)$ and $m_0P(S)/2$ do not intersect if $S$ is
  excessive. Thus Part~\eqref{i:ex} follows from Lemma~\ref{l:incl} as
  well.
\end{proof}

\section{Structures of Exceptional Sets}
\label{s:structure} We study of the structures of exceptional sets in
this section. As a result, we prove a rather curious fact: a spectral
class can only have either $1, 2$ or infinitely many elements. We
start by making the following conventions and definitions. Throughout
this section, $S$ denotes a balanced spectrum and $m$ denotes the
largest odd element of $S$, moreover:
\begin{itemize}
\item For a prime $p$, let $\gamma_p$ denote the largest integer such
  that $p^{\gamma_p} \in S_1$ and we write $\mu_p$ for $v_p(m)$. The
  {\bf excessive index of $p$ with respect to $S$} is defined to be
  $\epsilon_p:=\gamma_p -\mu_p$. Note that $\epsilon_p$ is always
  non-negative.

\item A prime $p$ is called an {\bf excessive prime} of $S$ if
  $\epsilon_p > 0$; otherwise $p$ is called a {\bf non-excessive
    prime} of $S$. We also say that $p$ is {\bf excessive
    (non-excessive) with respect to $S$} if it is an excessive (a
  non-excessive) prime of $S$. Note that $2$ is always a non-excessive
  prime of $S$. Also, every excessive prime of $S$ is in $S_1$ and
  every non-excessive prime of $S$ that is in $S_1$ divides $m$.

\item The {\bf excessive number} of $S$ is defined to be the product
  $e_S:=\prod p^{\epsilon_p}$ where $p$ runs through the primes. Note
  that $S$ is excessive if and only if $S$ has an excessive prime if
  and only if $e_S > 1$.

\item Every $a \in E(S)$ can be written as $e_an_a$ where $e_a$
  ($n_a$) is a product of (non-)excessive primes. The numbers $e_a$
  and $n_a$ are called the {\bf excessive part} and the {\bf
    non-excessive part} of $a$, respectively. Note that either $e_a$
  or $n_a$ can be $1$ but not both since $a > 1$.
\end{itemize}

The next two lemmas tell us what kind of factors that an element of
$E(S)$ can/must have.
\begin{lemma}
  \label{l:excessive} For every prime $p$, $p^{\epsilon_p}$ divides
  every element of $E(S)$.
\end{lemma}
\begin{proof} 
  For any $a \in E(S)$, since $p^{\gamma_p} \in S_1$, $p^{\gamma_p}
  \mid ma$ and so $p^{\epsilon_p} = p^{\gamma_p - \mu_p} \mid a$.
\end{proof}

\begin{lemma}
  \label{l:non-excessive} If $q$ is a non-excessive prime of $S$
  dividing some element of $E(S)$ then $q > p^{\mu_p}$ for any prime
  $p$ other than $q$.
\end{lemma}
\begin{proof} 
  Suppose $q$ is non-excessive and $q \mid a$ for some $a \in E(S)$.
  Then $q^{\mu_q+1}$ divides $ma$; moreover for any prime $p \neq q$,
  if $q < p^{\mu_p}$ then
  \[ 
    q^{\mu_q+1} \le \frac{mq}{p^{\mu_p}}< m < a.
  \] 
  But that means $q^{\mu_q+1} \in \fl{ma}{a}=S_1$, contradicting the
  fact that $q$ is non-excessive. Therefore, we must have $p^{\mu_p} <
  q$.
\end{proof}

\begin{proposition}
  \label{p:form} For every $a \in E(S)$, 
  \begin{enumerate}[(i)]
  \item \label{i:divide} $e_a$ is divisible by $e_S$.
  \item \label{i:either} $n_a$ is either $1$ or a power of the largest
    non-excessive prime of $S$ in $S_1$. In particular, if $n_a >1$
    then $S_1$ contains a non-excessive prime of $S$.
  \item \label{i:n_a>1} if $n_a > 1$, then $e_a =e_S$.
\end{enumerate}
\end{proposition}
\begin{proof} 
  By Lemma~\ref{l:excessive}, $e_S \mid a$. Since $e_S$ is a product
  of excessive primes, so in fact, $e_S$ divides $e_a$.

  By Lemma~\ref{l:facts}~\eqref{i:notprime} and~\eqref{i:proper},
  every prime factor of $a$ is in $S_1$. So $n_a$ is a product of
  non-excessive primes (of $S$) in $S_1$.  Therefore, $n_a=1$ if $S_1$
  contains no non-excessive primes of $S$. So suppose otherwise and
  let $q$ be the largest non-excessive prime of $S$ in $S_1$.  If
  $q_0$ is another non-excessive prime of $S$ dividing $n_a$, then
  $q_0 \in S_1$ and so $q_0 < q$. However, by
  Lemma~\ref{l:non-excessive} $q^{\mu_q} < q_0$ and this leads to a
  contradiction since $\mu_q \ge 1$. So we conclude that $n_a$ can
  have no prime factors other than $q$, therefore
  Part~\eqref{i:either} follows.

  Suppose $n_a >1$ then, by Part~\eqref{i:either}, $a$ is of the form
  $e_aq^{\beta}$ where $q$ is the largest non-excessive prime of $S$
  in $S_1$ and $\beta \ge 1$. For a prime $p$, let us write $\alpha_p$
  for $v_p(a)$. Consider the factor $p^{\mu_p + \alpha_p}$ of $ma$.
  For $p \neq q$, $p^{\alpha_p} \mid e_a$ and by
  Lemma~\ref{l:non-excessive} $p^{\mu_p} < q$. Therefore, $p^{\mu_p +
    \alpha_p} < a$ and so $p^{\mu_p + \alpha_p} \le m$ since $a \in
  E(S)$.  Consequently, $\mu_p + \alpha_p \le \gamma_p$, i.e.
  $\alpha_p \le \gamma_p - \mu_p = \epsilon_p$.  But by
  Part~\eqref{i:divide}, $\epsilon_p \le \alpha_p$. Therefore we
  conclude that $\alpha_p=\epsilon_p$ for every prime $p \neq q$ and
  hence $e_a=e_S$.
\end{proof} 
We should point out that it is possible for a balanced spectrum (other
than $\{1,2\}$) to have no non-excessive odd primes (see
Example~\ref{ex:no-ne}). Also, $n_a$ in the above proposition may
still be $1$ even $S_1$ contains a non-excessive prime of $S$.

Next we give a characterization of non-excessive spectra in terms of
exceptional sets.
\begin{theorem} 
  \label{th:E(S)non-excessive} A balanced spectrum $S$ is
  non-excessive if and only if $E(S)=\emptyset$ or
  $E(S)=\{q^{\mu_q+1}\}$ for some prime $q$. In particular, the size
  of the exceptional set of a non-excessive spectrum is at most one.
\end{theorem}
\begin{proof} 
  Suppose $S$ is non-excessive and $a \in E(S)$. Since every prime is
  non-excessive with respect to $S$, $e_a=1$ and so $n_a >1$.
  Therefore, by Proposition~\ref{p:form}~\eqref{i:either}, $a$ is a
  power of the largest non-excessive prime $q$ in $S_1$ (since $S$ is
  non-excessive, so in fact $q$ is outright the largest prime in
  $S_1$).  Since $a > m \ge q^{\mu_q}$, so on one hand $a \ge
  q^{\mu_q+1}$; on the other hand $q^{\mu_q+1} \mid ma$ but
  $q^{\mu_q+1} \notin S_1$ thus $a \le q^{\mu_q+1}$.  Therefore,
  $E(S)$ must be the singleton $\{q^{\mu_q+1}\}$ if it is non-empty.

  To show the other implication, let us note that if $E(S)$ is empty,
  then $S$ is non-excessive by
  Theorem~\ref{th:lsc-balanced}~\eqref{i:ex}. So let us assume
  $E(S)=\{q^{\mu_q+1}\}$ for some prime $q$. Then every $k \in S_1$
  divides $mq^{\mu_q+1}$ and $k < q^{\mu_q+1}$.  Thus $v_q(k) \le
  \mu_q$.  Moreover, for any prime $l$ other than $q$, $v_{l}(k) \le
  v_{l}(mq^{\mu_q+1}) = v_{l}(m)$.  Therefore, we conclude that $k
  \mid m$ and hence $S$ is non-excessive.
\end{proof}

The next result was a surprise to us.
\begin{theorem}
  \label{th:E(S)>1} $|E(S)| \le 2$.
\end{theorem}
\begin{proof} 
  Suppose $|E(S)| > 1$, then $S$ is excessive according to
  Theorem~\ref{th:E(S)non-excessive}. Let $p$ be an excessive prime of
  $S$. By Lemma~\ref{l:excessive}, $p^{\epsilon_p}$ divides every
  element of $E(S)$. Moreover, since $\epsilon_p \ge 1$, by
  Lemma~\ref{l:facts}~\eqref{i:proper} $p^{-1}E(S)$ and hence
  $p^{-\epsilon_p}E(S)$ is a subset of $S_1$.

  Let $(u_i)$ be the list of elements of $U:=p^{-\epsilon_p}E(S)$ in
  ascending order. For $i \ge 2$, since $m < p^{\epsilon_p}u_1 <
  p^{\epsilon_p}u_i \in E(S)$, therefore $p^{\epsilon_p}u_1$ must not
  divide $mp^{\epsilon_p}u_i$.  In other words, $u_1$ does not divide
  $mu_i$.  However, since $p^{\epsilon_p-1}u_1 \in p^{-1}E(S)
  \subseteq S_1$, $p^{\epsilon_p-1}u_1 \mid mp^{\epsilon_p}u_i$, i.e.
  $u_1 \mid mpu_i$. That means $v_p(u_1) \le v_p(mu_i)+1$ and for any
  prime $l$ other than $p$, $v_{l}(u_1) \le v_{l}(mpu_i) =
  v_{l}(mu_i)$. So the fact that $u_1$ does not divide $mu_i$ implies
  $v_p(u_1) > v_p(mu_i)$.  Therefore, we must have
  \[ 
  v_p(u_1) = v_p(mu_i) + 1 = \mu_p + v_p(u_i) + 1.
  \] 
  We claim that for $i \ge 2$, $u_i$ is not divisible by $p$. If not,
  then $v_p(u_1) > \mu_p+1$ and so $p^{\gamma_p+1}$ would be a proper
  factor of $p^{\epsilon_p}u_1$ and hence, by
  Lemma~\ref{l:facts}~\eqref{i:proper}, an element of $S_1$, a
  contradiction. Therefore, we conclude that
  \[ 
  E(S) = p^{\epsilon_p}\{p^{\mu_p+1}v_1, u_2, \ldots, u_t\} =
  \{p^{\gamma_p+1}v_1, p^{\epsilon_p}u_2, \ldots, p^{\epsilon_p}u_t\}
  \] 
  for some $v_1, u_2, \ldots, u_t$ not divisible by $p$. Since
  $p^{\gamma_p+1} \notin S$, then again by
  Lemma~\ref{l:facts}~\eqref{i:proper} $v_1=1$ and so $p^{\gamma_p+1}$
  is the least element of $E(S)$.

  As elements of $S_1$ the $u_i$'s all divide $mp^{\gamma_p+1}$. But
  for $i \ge 2$, $u_i$ and $p$ are relatively prime so $u_i$ must
  divide $m$. Therefore, if $|E(S)| > 2$, then we would have
  $p^{\epsilon_p}u_2 \mid mp^{\epsilon_p}u_3$ and $m <
  p^{\epsilon_p}u_2 < p^{\epsilon_p}u_3$, contradicting
  $p^{\epsilon_p}u_3 \in E(S)$. So we conclude that $|E(S)| \le 2$.
\end{proof}

\begin{theorem}
  \label{th:size} Every spectral class has either $1,2$ or infinitely
  many elements. Moreover, for any $n$,
  \begin{itemize}
  \item $|L(n)|=1$ if and only if $\lspec(n)$ is lopsided or excessive
    with an exceptional set of size $1$.

  \item $|L(n)|=2$ if and only if $\lspec(n)$ is excessive with an
    exceptional set of size $2$.

  \item $L(n)$ is infinite if and only if $\lspec(n)$ is unmixed or
    non-excessive.
  \end{itemize}
\end{theorem}
\begin{proof} 
  Suppose a spectral class $L(n)$ is finite, then $S:=\lspec(n)$ must
  be either lopsided or excessive (Theorem~\ref{th:lsc-unmixed}
  and~\ref{th:lsc-balanced}~\eqref{i:non-ex}). In the former case,
  $|L(n)|=|\{n\}|=1$ (Theorem~\ref{th:lsc-lopsided}); in the latter
  case, $1 \le |L(n)|=|E(S)| \le 2$ according to
  Theorem~\ref{th:lsc-balanced}~\eqref{i:ex} and
  Theorem~\ref{th:E(S)>1}. So we establish the first statement.  The
  three equivalences are simply reorganizing what we have already
  proved in Theorem~\ref{th:lsc-lopsided}, \ref{th:lsc-unmixed}
  and~\ref{th:lsc-balanced}.
\end{proof}

We conclude this section with a few more precise descriptions of the
exceptional sets.

\begin{proposition}
  \label{p:s1} If $e_S > m$, then $E(S)=\{e_S\}$.
\end{proposition}
\begin{proof} 
  If $e_S > m$, then in particular $e_S$ is greater than $1$ so $S$ is
  excessive and hence $|E(S)| \ge 1$. Let $a \in E(S)$, by
  Lemma~\ref{l:excessive}, $e_S \mid a$. However, since $e_S > m$,
  therefore by Lemma~\ref{l:facts}~\eqref{i:proper}, $e_S = a$.
\end{proof}
\begin{proposition}
  \label{p:E(S)=1} Suppose $|E(S)|=1$ then the unique element of
  $E(S)$ is either a product of excessive primes or of the form
  $e_Sq^{\beta}$ {\upshape(}$\beta \ge 1${\upshape)} where $q$ is the
  largest non-excessive primes of $S$ in $S_1$.
\end{proposition}
\begin{proof}
  Let $a$ be the unique element of $E(S)$. If $a$ is not a product of
  excessive primes then by Proposition~\ref{p:form}~\eqref{i:proper}
  and~\eqref{i:n_a>1}, $a$ is of the form $e_Sq^{\beta}$ with $\beta
  \ge 1$.
\end{proof}

\begin{proposition}
  \label{p:s2} Suppose $|E(S)| > 1$ then $S$ has a unique excessive
  prime $p$ and $E(S)$ is of the form $\{p^{\gamma_p+1},
  p^{\epsilon_p}q^{\beta}\}$ {\upshape(}$\beta \ge 1${\upshape)} where
  $p$ and $q$ are the two largest primes in $S_1$. Moreover,
  \begin{enumerate}[(i)]
  \item \label{i:p>q} if $p$ is the largest prime in $S_1$, then
    $\epsilon_p = \gamma_p$.
  \item \label{i:q>p} if $p$ is not the largest prime in $S_1$, then
    $\beta=1$.
  \end{enumerate}
\end{proposition}
\begin{proof} 
  We have already proved (Theorem~\ref{th:E(S)>1}) that $|E(S)| > 1$
  implies $E(S)$ is of the form $\{p^{\gamma_p+1}, p^{\epsilon_p}u\}$
  where $p$ is an excessive prime of $S$, $u > p^{\mu_p+1}$ and $u$ is
  not divisible by $p$. By Lemma~\ref{l:excessive}, every excessive
  prime of $S$ divides $p^{\gamma_p+1}$, therefore $p$ is the only
  excessive prime of $S$ and so $u$ is the non-excessive part of
  $p^{\epsilon_p}u$.  Since $u > 1$, by Proposition~\ref{p:form} $u$
  is a positive power of $q$ where $q$ is the largest non-excessive
  prime of $S$ in $S_1$.  Therefore, we conclude that $E(S)$ must be
  of the form $\{p^{\gamma_p+1}, p^{\epsilon_p}q^{\beta}\}$ for some
  $\beta \ge 1$.

  Since $q$ is the largest non-excessive prime in $S_1$ and $p$ is the
  unique excessive prime in $S_1$, therefore if $p$ is the largest
  prime in $S_1$ then $q$ must be the second largest prime in $S_1$.
  Since $q$ divides $p^{\epsilon_p}q^{\beta} \in E(S)$, by
  Lemma~\ref{l:non-excessive}, $q > p^{\mu_p}$. Therefore, $\mu_p$
  must be $0$, i.e. $\epsilon_p=\gamma_p$. This completes the proof
  of~\eqref{i:p>q}.

  Suppose $p$ is not the largest prime in $S_1$, then the largest
  prime in $S_1$ is non-excessive and so must be $q$. We claim that in
  this case $\beta$ is actually $1$. First, note that $m <
  p^{\gamma_p+1} < p^{\gamma_p}q$ and since $p^{\gamma_p}q$ divides
  $mp^{\epsilon_p}q^{\beta}$, $p^{\epsilon_p}q^{\beta} \le
  p^{\gamma_p}q$ i.e.  $q^{\beta-1} \le p^{\mu_p}$. But by
  Lemma~\ref{l:non-excessive}, we also have $p^{\mu_p} < q$.
  Therefore, $\beta$ must be $1$.

  To finish the proof of~\eqref{i:q>p}, we argue that $p$ must be the
  second largest prime in $S_1$. Since $q^{\mu_q+1} \mid
  mp^{\epsilon_p}q$ and $q$ is non-excessive, $q^{\mu_q +1} \ge
  p^{\epsilon_p}q$. In other words, $p^{\mu_p}q^{\mu_q} \ge
  p^{\gamma_p}$. So if there were a prime $l \in S_1$ strictly in
  between $p$ and $q$ then $p^{\mu_p}lq^{\mu_q} > p^{\gamma_p+1} > m$.
  But since $l$ is non-excessive, it divides $m$ and therefore
  $p^{\mu_p}lq^{\mu_q} \mid m$, a contradiction.
\end{proof}

\section{Examples and Algorithms} 
\label{s:e-a} We give some examples here to illustrate the results in
previous sections.
\begin{example} 
  The smallest number with an unmixed spectrum is $1$. Since the
  spectrum of $1$ is $\{1\}$, it follows from
  Theorem~\ref{th:lsc-unmixed} that $L(1)$ is the set of powers of
  $2$.
\end{example}

\begin{example} 
  The smallest number with a balanced spectrum is $3$. The spectrum of
  $3$ is $\{1,2\}$ which is non-excessive and has an empty exceptional
  set. By Theorem~\ref{th:lsc-balanced}~\eqref{i:non-ex}, $L(3)$ is
  the set of odd primes.
\end{example}

\begin{example} 
  The smallest number with a lopsided spectrum is $9$. The spectrum of
  $9$ is $\{1,2,3\}$. By Theorem~\ref{th:lsc-lopsided}, $L(9) =\{9\}$.
\end{example}

\begin{example} 
  The number $21$ is the smallest number with a spectrum that has a
  non-empty exceptional set. The spectrum of $21$ is $\{1,2,3,6\}$.
  It is non-excessive with exceptional set $\{9\}$. By
  Theorem~\ref{th:lsc-balanced}~\eqref{i:non-ex},
  \[ L(21) = \{27\} \cup \{3p: p \text{ prime} > 6\}. \]
\end{example}

\begin{example} 
  The smallest number with an excessive spectrum is $75$. The
  exceptional set of $\lspec(75)=\{1,2,3,5,6,10\}$ is $\{15\}$. By
  Theorem~\ref{th:lsc-balanced}~\eqref{i:ex},
  $L(75)=\{75\}$.
\end{example}

\begin{example} 
  The smallest number with two elements in the exceptional set of its
  spectrum is $175$. The spectrum of $175$ is $\{1,2,5,7,10,14\}$. The
  exceptional set of $\lspec(175)$ is $\{25,35\}$ (c.f.
  Proposition~\ref{p:s2}~\eqref{i:q>p}). By
  Theorem~\ref{th:lsc-balanced}~\eqref{i:ex},
  \[ L(175)=\frac{14}{2}\{25,35\} = \{175,245\}. \]
\end{example}

\begin{example} 
  The proof of Theorem~\ref{th:E(S)>1} would be considerably simpler
  if every excessive spectrum contains an odd prime not dividing its
  largest odd element. However, it is not always the case. The
  smallest number with a spectrum witnessing this fact is $2673=3^5
  \cdot 11$.  The set of odd elements of the spectrum of $2673$ is
  $\{1,3,9,11,27,33\}$. It contains two primes, $3$ and $11$, both of
  them divide $33$.
\end{example}

\begin{example}
  \label{ex:no-ne} The number $9261=3^3\cdot 7^3$ is the smallest
  number such that its spectrum contains an odd prime and every odd
  prime in its spectrum is excessive. The set of odd elements of
  $\lspec(9261)$ is $\{1,3,7,9,21,27,49,63\}$. It contains two primes,
  $3$ and $7$, both of them have excessive index $1$. Let us explain
  how we found this example.  Suppose $S=\lspec(n)$ has the required
  property. Then by Proposition~\ref{p:s2}, $E(S)$ must be a
  singleton. Also, $m:=\max S_1$ cannot be a prime power, otherwise
  the prime of which $m$ is a power will be non-excessive.  So $m$ has
  at least two prime factors, say $p_1 < p_2$. Since $p_2$ is
  excessive, therefore $p_2^2 \in S_1$. But that means $m > p_2^2 >
  p_1p_2$ and so $m=p_1p_2c$ for some odd number $c > 1$. By
  Lemma~\ref{l:excessive}, the unique element $a \in E(S)$ is of the
  form $p_1p_2d$. Since it must be greater than $2m$, $d > 2c$. At
  this juncture, we make a guess: suppose $m=p_1^2p_2$ and
  $a=p_1p_2^2$.  Then $n=(p_1p_2)^3$ and we need $p_2/p_1 = a/m > 2$.
  Minimizing $n$ subjected to the inequality yields $p_1=3, p_2=7$ and
  so $n=(21)^3 = 9261$.

  Now let $n_0$ be the smallest number such that its spectrum contains
  an odd prime and every odd prime in its spectrum is excessive. The
  argument above shows that $n_0 \le 9261$ and is the product of two
  numbers of the forms $p_1p_2c$ and $p_1p_2d$ with $d > 2c$. Since
  $2(p_1p_2c)^2=2m^2 < n_0 \le 9261$, $c^2 < 9261/2(15)^2$.
  Therefore, $c$ must be $3$ and $d \ge 7$. From this, we see that
  $p_1^2p_2^2 \le 9261/21= (21)^2$. If $p_1=3,p_2=5$, then $7 \le d
  \le 9261/(3\cdot (15)^2)$, i.e.  $d=7,9,11$ or $13$. But none of
  these choices produces a balanced spectrum.  So $p_1=3,p_2=7$, this
  forces $d=7$ and hence $n_0=9261$.
\end{example}

\begin{example} 
  \label{ex:beta>1} The exponent $\beta$ in Proposition~\ref{p:s2} can
  be greater than $1$. Let us find the smallest witness, say $n_0$, of
  this fact. First, let $n$ be a number with $S=\lspec(n)$ witnessing
  $\beta > 1$. Let $p$ be the unique excessive prime of $S$ and $q$ be
  the largest non-excessive prime in $S_1$.  For simplicity, let us
  write $\gamma$ for $\gamma_p$. By
  Proposition~\ref{p:s2}~\eqref{i:p>q}, $p >q$, $E(S)=\{p^{\gamma+1},
  p^{\gamma}q^{\beta}\}$ and $p$ does not divide $m:=\max S_1$.
  Therefore, $m$ is of the form $cq^{\beta+\kappa}$ where $\kappa \ge
  0$ and $c$ is not divisible by either $p$ or $q$. The number of
  factors of the two elements in $mE(S)$ are both equal to the size of
  $S$; by equating them, we get $\beta\gamma = \kappa+1$. Clearly, one
  solution of this equation with $\beta >1$ is $\beta =2,
  \gamma=\kappa=1$. With these values, $E(S)=\{p^2, pq^2\}$. So the
  choices of $p$ and $q$ have to meet the following inequalities:
  \begin{equation}
    \label{eq:mincond}
    (q <) p< q^2\ \text{and}\ (m = )cq^3 < p^2/2.
  \end{equation}
  Minimizing $n=cq^3p^2$ subjected to these inequalities yields $c=1$,
  $p=17$ and $q=5$. So $n=5^3\cdot 17^2=36125$ and the exceptional set
  of its spectrum is $\{17^2, 17\cdot 5^2\}$. It is indeed an example
  with $\beta(=2) >1$.

  To show that $n_0$ is $36125$, let us note that, from the argument
  above, $n_0$ is of the form $cq^{\beta+\kappa}p^{\gamma+1}$ with
  $cq^{\beta+\kappa} < p^{\gamma+1}/2$. Therefore, $n_0 >
  2c^2q^{2(\beta+\kappa)}$. Since $\beta\gamma=\kappa+1$ and $\beta
  \ge 2$, $36125 \ge n_0 > 2q^6$.  That means $q \le 5$. But one
  quickly rules out the possibility that $q=3$, since no choice of $p$
  would satisfy the inequalities in~\eqref{eq:mincond}. Hence $q=5$,
  and it follows that the minimal choices for $c$ and $q$ are $1$ and
  $17$, respectively. This completes the proof.

  To find an exceptional set with both $\gamma$ and $\beta > 1$ will
  lead one to the number $21434375=5^5\cdot 19^3$. The exceptional set
  of its spectrum is $\{19^3, 19^2\cdot 5^2\}=\{6859,9025\}$. An
  analysis similar to the one given above shows that $21434375$ is
  indeed the smallest possible choice. We will leave the verification
  to the reader this time.
\end{example}

Finally, we give two algorithms answering the questions that we posed
in the introduction. Algorithm~\ref{al:spectral-class} computes from a
given number $n$ its spectral class $L(n)$.  Its correctness is
guaranteed by Theorem~\ref{th:lsc-lopsided},~\ref{th:lsc-unmixed}
and~\ref{th:lsc-balanced}. Algorithm~\ref{al:lspec} computes, using
Algorithm~\ref{al:spectral-class}, from a finite set of numbers $S$
the set of numbers with spectrum $S$. In particular,
Algorithm~\ref{al:lspec} returns the empty set if $S$ is not a
spectrum. Here is the strategy: computes from $S$ a number $n$ such
that $S=\lspec(n)$ if $S$ is a spectrum.  The algorithm then returns
either $L(n)$ or the empty set depending on whether the equality holds
or not.  We have implemented both algorithms using PARI/GP script.

\begin{algorithm}
  \caption{Compute the spectral class of a natural number}
  \label{al:spectral-class}
  \begin{algorithmic}
    \REQUIRE A natural number $n$
    \ENSURE  $L(n)$ the spectral class of $n$
    \STATE $S:=\lspec(n)$; $S_i:=\{a \in S \colon a \equiv i
      \pmod{2}\}$ ($i=0,1$).
    \IF{$|S_0|=0$}  
      \STATE $m_1:=\max S_1$; $\nu:=$ the least integer such that
      $2^{\nu+1} > m_1$.
      \RETURN $\{2^{\nu+i}m_1 \colon i \ge 0\}$.
      \ELSIF{$|S_0| < |S_1|$}
         \RETURN $\{n\}$.
         \ELSIF{$S$ is non-excessive}
            \STATE $m_0:= \max S_0$.
           \RETURN $\frac{1}{2}m_0(P(S) \cup E(S))$.
         \ELSE
            \STATE $m_0:= \max S_0$.
            \RETURN $\frac{1}{2}m_0E(S)$.
    \ENDIF
  \end{algorithmic}
\end{algorithm}

\begin{algorithm}
\caption{Compute the set of numbers with a given set as
  spectrum}
\label{al:lspec}
\begin{algorithmic}
  \REQUIRE A finite set of natural numbers $S$.
  \ENSURE The set of natural numbers with length spectrum $S$.
  \STATE $S_i :=\{a \in S \colon a \equiv i \pmod{2}\}$ ($i=0,1$).
  \IF{$|S_0|= 0$}
     \STATE $m_1:= \max S_1$; 
     $\nu:=$ the least integer such that $2^{\nu+1} > m_1$; 
     $n:=2^{\nu}m_1$.
     \ELSIF{$|S_0| < |S_1|$}
        \STATE $n:=\frac{1}{2}\min S_0\min D(S)\max D(S)$.
        \ELSIF{$|E(S)| \neq 0$}
        \STATE $m_0:=\max S_0$; $n:=\frac{1}{2}m_0\min E(S)$.
        \ELSE
        \STATE $m_0:=\max S_0$; 
        $p:=$ the least prime $> m_0$; $n:=\frac{1}{2}m_0p$.
  \ENDIF
  \IF{$\lspec(n)=S$} 
    \RETURN $L(n)$.
  \ELSE 
  \RETURN $\emptyset$.
  \ENDIF
\end{algorithmic}
\end{algorithm}

\end{document}